\documentclass[aap]{imsart}

\usepackage{amsthm,amsmath,amssymb,natbib}

\arxiv{math.PR/0000000-***}

\startlocaldefs

\numberwithin{equation}{section}

\newtheorem{theorem}{Theorem}[section]
\newtheorem{lemma}[theorem]{Lemma}
\newtheorem{proposition}[theorem]{Proposition}

\newtheorem{remark}[theorem]{Remark}

\newcommand{\bpf}[1][Proof]{{\noindent {\sc #1: }}}
\newcommand{\epf}{{\hfill $\square$}\vspace{.7cm}}

\newcommand{\abs}[1]{\lvert#1\rvert}
\newcommand{\qvar}[1]{\langle#1\rangle}

\newcommand{\Dql}{D_{\mathrm{QL}}}

\newcommand{\eg}{{\textrm{e.g.\ }}}
\newcommand{\ie}{{\textrm{i.e.\ }}}
\newcommand{\viz}{{\textrm{viz.\ }}}
\newcommand{\rmd}{{\mathrm{d}}}
\newcommand{\rme}{{\mathrm{e}}}
\newcommand{\rmi}{{\mathrm{i}}}

\newcommand{\EE}{{\mathbb E}}
\newcommand{\NN}{{\mathbb N}}
\newcommand{\PP}{{\mathbb P}}

\newcommand{\RR}{{\mathbb R}}
\newcommand{\TT}{{\mathbb T}}
\newcommand{\ZZ}{{\mathbb Z}}

\newcommand{\cB}{{\mathcal B}}

\newcommand{\cF}{{\mathcal F}}
\newcommand{\cM}{{\mathcal M}}
\newcommand{\cN}{{\mathcal N}}

\newcommand{\cT}{{\mathcal T}}

\newcommand{\mfm}{{\mathfrak{m}}}

\endlocaldefs

\begin{document}

\begin{frontmatter}

\title{Diffusion limit for many particles
       in a periodic stochastic acceleration field}
\runtitle{Diffusion limit for stochastic acceleration}

\begin{aug}


\author{\fnms{Yves} \snm{Elskens}\corref{}\ead[label=e1]{yves.elskens@univ-provence.fr}} \and
\author{\fnms{Etienne} \snm{Pardoux}\ead[label=e2]{pardoux@cmi.univ-mrs.fr}}

 \affiliation{CNRS--universit\'e de Provence}

 \address{Physique des interactions ioniques et mol\'eculaires,
          UMR 6633, \\
          Equipe turbulence plasma,
          case 321 campus Saint--J{\'e}r{\^o}me, 
          av.~escadrille Normandie--Niemen, \\ FR-13397 Marseille cedex 13 \\
          \printead{e1}}

 \address{Laboratoire d'analyse, topologie et probabilit\'es,
          UMR 6632, \\
          centre de math\'ematique et d'informatique, \\
          rue F.~Joliot--Curie, 39,
          FR-13453 Marseille cedex 13 \\
          \printead{e2}}

\runauthor{Y.~Elskens and E.~Pardoux}
\end{aug}

\begin{abstract}
  The one-dimensional motion of any number $\cN$ of particles
  in the field of many independent waves (with strong spatial correlation)
  is formulated as a second-order system of stochastic differential equations,
  driven by two Wiener processes.
  In the limit of vanishing particle mass $\mfm \to 0$,
  or equivalently of large noise intensity,
  we show that the momenta of all $\cN$ particles converge weakly to $\cN$ independent
  Brownian motions, and this convergence holds even if the noise is periodic.
  This justifies the usual application of the diffusion equation to a family of
  particles in a unique stochastic force field.
  The proof rests on the ergodic properties of the relative velocity of two particles
  in the scaling limit.
\end{abstract}






\end{frontmatter}




\section{Introduction}

The motion of a particle in the field of many waves \cite{DG03,DoveilMacor06,Tsunoda91}
is a fundamental process in classical physics, the understanding of which is a
prerequisite to the analysis of many plasma and fluid phenomena. In one space dimension,
it can be described by the Hamiltonian model
\begin{equation}
  H
  =
  \frac{p^2}{2{\mfm}}
  + \sum_{m=1}^\cM A_m \cos (k_m q - \omega_m t - \varphi_m)
  \label{Hns}
\end{equation}
where the particle with mass ${\mfm}$ has position $q$ and momentum $p$, while the force
field derives from a potential with time Fourier components $A_m \rme^{\rmi \varphi_m}$.
The wave field comprises $\cM$ waves, with a smooth dispersion relation associating a
wavenumber $k_m$, a pulsation $\omega_m$ and a phase velocity $v_m = \omega_m / k_m$ to
each wave -- usually determined by fixed properties of the environment, such as the
geometry of the domain where waves propagate (then wavenumbers $k_m$ and pulsations
$\omega_m$ are discrete). The complex amplitudes $A_m \rme^{\rmi \varphi_m}$ are more
easily tuned by the experimenter or affected by simple changes in the environment.

The dynamical systems approach to this problem discusses the particle motion after
prescribing a single choice for each wave complex amplitude. As it would be quite
exceptional to control all waves (though this is \eg the assumption underlying the
standard map, see \cite{BE98Std} for a discussion), physicists often turn to a
probabilistic description of the dynamics, considering an ``ensemble'' of realizations
$(A_m,\varphi_m)$. Various arguments are then invoked to reduce the particle evolution
equations to a stochastic differential equation, often driven by a ``white noise''. This
results in somewhat tractable models (see \eg \cite{Krommes} about the validity of such
derivations).

In this paper we focus on two issues. First, a random field characterized by
$(A_m,\varphi_m)$ for a given dispersion relation with discrete frequency spectrum may
be periodic in time~: may the force on the particle be considered as independent over
several time periods in a genuine limit~? Second, may one consider several particles
subject to the same wavefield as independent in a genuine limit~? The latter issue
underlies the frequent application of the Fokker--Planck equation to the evolution of a
family of particles in a single turbulent wavefield -- though a priori one can only
grant that the diffusion equation describes the evolution of the distribution of a
single particle for an ensemble of wavefield samples.

From a more general perspective, this work also relates to the issue of ``propagation of
chaos'' in statistical physics \cite{Kac56,Kac59}, an aspect of Hilbert's 6th problem~:
how does chaotic dynamics enable a system, in which initial data are independent
(``random'') but the evolution may generate correlations, to behave as if the evolution
regenerated independence (``randomness'') or destroyed correlations~? Here, how do two
Wiener processes, fully describing a prescribed ``turbulent'' environment, generate
$\cN$ independent Brownian motions for particles~?

A further motivation for the present work is that physics literature most often focuses
on the evolution of particle distribution functions, \eg by showing that they obey a
Fokker--Planck equation, and on instantaneous observables such as $p(t_1)$ for given
$t_1$ (pointwise in time). However, the notion of a diffusion process implies rather a
measure on the set of trajectories, \viz functions $p(\cdot)$ (globally in time). Here
we shall show how our model implies that an arbitrary number of trajectories in a single
realization of the dynamics do, jointly, admit the Wiener measure description.

In sections \ref{back} and \ref{mainres} we motivate the mathematical model more
precisely and state our main results, which are proved in the subsequent sections. The
crucial Theorem \ref{theth} is an ergodic theorem for a rescaled process, implying that
a process $V_t$, describing the relative velocity of one particle with respect to
another one (or to its own earlier motion), converges weakly to a Brownian motion. The
difficulty in proving the ergodic theorem is that the invariant measure of the diffusion
process is infinite (it is the Lebesgue measure $\rmd x \rmd y$), and we must estimate a
continuous additive functional generated by a function (namely $f(x,y) = \sin^2 x$)
which is not integrable with respect to this measure (it is only locally integrable).
This weak convergence then implies the final many--particle result (Theorem
\ref{NpartGauss}) by a straightforward application of the L\'evy characterization of
Brownian motion. Section \ref{persp} outlines implications and possible extensions to
this work.

\section{Physical background}
\label{back}

Because the waves have different frequencies and velocities, it is generally unrealistic
to assume their phases to be correlated. Their intensities are more easily observed, but
both in nature and in the laboratory the accumulation of statistical data on waves often
involves only their average power spectra, not the detailed intensity data for each
measurement run. We assume here that these complex amplitudes are random data, and
investigate the statistics of the particle motion in the resulting time--dependent
random field. This dynamics is a ``stochastic acceleration problem'' for a ``passive
particle'' in weak plasma turbulence
\cite{Dimits86,DG82,Maasjost,Sturrock66,VdEijnden97}, and its understanding is a
prerequisite to a proper analysis of the case where the particle motion feeds back on
the wave evolution \cite{DoxasCary,EEbook}.

The Hamiltonian \eqref{Hns} generates equations of motion
\begin{eqnarray}
  \dot q & = & p / {\mfm}
  \label{dotx} \, ,
  \\
  \dot p & = & \sum_m k_m A_m \sin (k_m q - \omega_m t - \varphi_m)
  \label{dotv}  \, .
\end{eqnarray}

An important observation \cite{Chirikov79,Escande85} on the motion of a particle in the
field (\ref{dotv}) is locality in velocity~: the evolution of the particle when it has
velocity $\dot q = v$ depends only weakly on the waves with a Doppler-shifted frequency
$\omega_m - k_m v$ much larger than their trapping oscillation frequency $k_m
\sqrt{A_m/{\mfm}}$. In particular, for a two-wave system the resonance overlap parameter
\begin{equation*}
  s_{1,2}
  =
  \frac {2 \sqrt{A_1 / {\mfm}} + 2 \sqrt{A_2 / {\mfm}}}
        {\abs{\omega_2 / k_2 - \omega_1 / k_1}}
\end{equation*}
becomes unity when there exists a velocity $u = \omega_2/k_2 - 2 \sqrt{A_2/{\mfm}}=
\omega_1/k_1 + 2 \sqrt{A_1/{\mfm}}$ (with $k_1>0$, $k_2 > 0$, $\omega_2 > \omega_1$).
For many waves with overlap parameters $s \gg 1$, the relevant phase velocity range for
waves influencing the particle is a ``resonance box'', with a width scaling as
$(A/{\mfm})^{2/3}$. \cite{BE97,BE98}

A good approximation to typical wave dispersion relations in the strong overlap limit,
after a Galileo change of reference frame (see \eg sec.~6.7 in \cite{EEbook}), is
\begin{equation}
  k_m = k_0 , \ \omega_m = 2 \pi (m-\cM/2) / \cT  ,
  \label{BEspectrum}
\end{equation}
for some $k_0$, $\cT$. Then, in the limit $\cM \to \infty$, the equations of motion
yield for $A_m \rme^{\rmi \varphi_m} = A_0$ (with real $A_0$) the well-known standard
map \cite{BE98Std}. The case where phases $\varphi_m$ are independent random variables
uniformly distributed on the circle $[0, 2 \pi]$, while $A_m = A_0$ is given, was
investigated notably by Cary, Escande, Verga and B\'enisti \cite{BE97,CEV,EEbook} and
occurs in the context of the random phase approximation.

To the extent that the phases and amplitudes of the waves are independent random
variables, the physicist usually views the force (\ref{dotv}) as a mollification of a
white noise, with amplitude $\sigma = \sqrt{\EE (k_m^2 A_m^2)} = |k_m| \sqrt{\EE
A_m^2}$, where the relevant mode $m$ is the one nearest to the current particle velocity
(the mathematical expectation $\EE$ is called ensemble average with respect to wave
amplitudes and phases).\footnote{Phases do not appear in $\sigma$ (nor in $s$) because
the relative phase of two waves $\varphi_m + \omega_m t - \varphi_n - \omega_n t$ varies
uniformly over time (hence $\varphi_m - \varphi_n$ can be absorbed in the choice of the
time origin).}

This is the core of quasilinear theory \cite{Drummond62,Pesme94,Romanov61,Vedenov62}.
With some care, one interprets (\ref{dotx})--(\ref{dotv}) as a stochastic differential
equation~; this applies in the case $m \in \ZZ$ for dispersion relation
(\ref{BEspectrum}) with Gaussian independent complex amplitudes such that $\EE A_m^2 =
k_0^{-2} \sigma^2$. The particle velocity then has a Brownian evolution, so that for $t,
t' \in [0,\cT]$
\begin{equation}
  \EE (p_t - p_0)(p_{t'} - p_0) = \Dql \min(t,t')
  \label{v2tql}
\end{equation}
with the quasilinear diffusion coefficient $\Dql = \sigma^2 \cT/2$. However, the
particle evolution for $t > \cT$ may show a strong correlation to its motion for $0 \leq
t \leq \cT$ because the waves are periodic in time \cite{E07,E08}, and the dispersion
relation (\ref{BEspectrum}) may generate a strong spatial correlation between the
motions of two particles because all waves acting on a particle at any time have the
same wavelength.\footnote{~There is a large body of literature on the case of incoherent
waves with no dispersion relation. Then the sum $\sum_m$ becomes a double sum
$\sum_{m,n}$ and one varies wavenumbers $k_n$ independently from pulsations $\omega_m$.
This space-time stochastic environment is more noisy than our model and may also be
considered to motivate a quasilinear approximation.}

Set $k_0 = 1$ and $\cT = 2 \pi$ by the choice of space and time units. The particle
phase space is $\TT \times \RR$, where $\TT$ is the circle modulo $2\pi$, and the
particle equation of motion reads, with initial data $(q(0),\dot q(0)) = (q_0, \dot
q_0)$,
\begin{equation}
  \begin{aligned}
    \ddot q &=  {\mfm}^{-1} \sum_{m = -\mu}^{\mu'-1} A_m \sin (q - m t + \varphi_m)
    \\
            &= {\mfm}^{-1} \sum_m A_m \cos (m t - \varphi_m) \sin q
             - {\mfm}^{-1} \sum_m A_m \sin (m t - \varphi_m) \cos q
  \end{aligned}
  \label{eqq2}
\end{equation}
in a galilean frame moving at a velocity inside the spectrum of wave phase velocities
($\mu + \mu' = \cM$). We assume $\mu \gg 1$, $\mu' \gg 1$. For $A_m = A_0$ given and
independent random phases (uniform on the circle), in the limit $A_0/{\mfm} \to \infty$,
with $({\mfm}/A_0)^{2/3} \cM \gtrsim 10$, B\'enisti and Escande \cite{BE97,BE98} have
shown that the particle momentum $p = {\mfm}\dot q$ follows essentially a Brownian
motion, with diffusion constant given by the quasilinear estimate
\begin{equation*}
  D_{\mathrm  QL} = \pi A_0^2
\end{equation*}
as long as the motion does not approach the boundaries of the wave velocity spectrum.
Then the particle momenta $p$ for an ensemble of independent realizations of the system
will be described by a distribution function verifying the Fokker--Planck equation
\begin{equation}
  \partial_t f = {D_{\mathrm QL} \over 2} \partial_p^2 f
  \label{FokPla} \,
\end{equation}
even for $t > \cT$. Numerical simulations \cite{E08} show similar behaviour for i.i.d.\
complex Gaussian random variables $A_m \rme^{\rmi \varphi_m}$. The present work
establishes a rigorous version of this result in the frame of stochastic processes.

\section{Main results}
\label{mainres}

We first let $\min(\mu,\mu')\to\infty$ in the model, taking $A_m \cos \varphi_m$, $A_m
\sin \varphi_m$ as i.i.d.\ Gaussian random variables with zero expectation and $\EE(A_m
\cos \varphi_m)^2 = \EE(A_m \sin \varphi_m)^2 = 1/2$ for all $m \in \ZZ$. In particular
this implies that phases $\varphi_m$ are i.i.d.\ uniformly on $\TT$. To follow earlier
practice, we now set ${\mfm} = 1/A$. Then formally (\ref{eqq2}) becomes the Stratonovich
stochastic differential equation for $0 \le t \le 2 \pi$
\begin{eqnarray}
  \rmd Q_t &=& A P_t \ \rmd t \, ,
  \label{eqq2b}
  \\
  \rmd P_t &=& \sin(Q_t) \circ \rmd C_t + \cos(Q_t) \circ \rmd S_t \, ,
  \label{eqp2b}
\end{eqnarray}
with initial data $Q_0 = q_0$, $P_0 = p_0 = {\mfm} {\dot q_0} = {\dot q_0}/A$~; here,
from \cite{Kahane}, $\pi^{-1/2}(C,S)$ is a standard 2-dimensional Brownian motion. In
other words, $C$ and $S$ are martingales, with $C_0 = S_0 = 0$ and
\begin{equation}
  \qvar{C}_t = \qvar{S}_t = \pi t  \, , \
  \qvar{C,S}_t = 0 \, .
  \label{EECS}
\end{equation}
Note that the Stratonovich and It\^o integrals define the same solutions for this
system, and that the vector fields $(\sin q)\partial_p$ and $(\cos q)\partial_p$
commute.

For $0 \leq t \leq 2 \pi$, it is clear that $P$ is a Brownian motion for any value of
$A\geq 0$ (see the first lines of Sec.~\ref{proofNpart}). However, the model
\eqref{eqq2} defines a dynamical system for all times $0 \leq t < \infty$, and one may
wonder how the initial stochastic behaviour over $[0,2\pi]$ extends for longer times.
Formally, one solves then \eqref{eqq2b}--\eqref{eqp2b} with the periodized field, \ie
with the continuous processes defined by
\begin{eqnarray}
  \rmd C_{t+2k\pi} = \rmd C_t \, & , &
  \rmd S_{t+2k\pi} = \rmd S_t
  \label{perCS}
\end{eqnarray}
for $k \in \ZZ$. In other words, $C_t - \frac{t}{2\pi}C_{2\pi}$ and $S_t -
\frac{t}{2\pi}S_{2\pi}$ are independent Brownian bridges repeated periodically for $t
\in \RR$, while $C_{2\pi}$ and $S_{2\pi}$ are independent Gaussian random variables with
expectation 0 and variance $2\pi^2$.

For this extended process, the wave field acting on the particle for $t \not \in
[0,2\pi]$ is not stochastically independent from the wave field acting during
$[0,2\pi]$. Therefore one does not expect the particle momentum to proceed as a Brownian
motion for all times, and indeed for $A$ small enough the velocity $\dot q$ may remain
bounded in a narrow interval for all times. This is easily seen numerically and can be
attributed to the existence of Kolmogorov-Arnol'd-Moser invariant tori in the
3-dimensional extended phase space with coordinates $(t,q,\dot q)$.

On the other hand, for large $A$, the dynamics viewpoint \cite{BE97,BE98,EEbook}
suggests that the nonlinearity in the equations of motion (due to trigonometric
functions of $Q$) may enable a decorrelation of the force over the period $\cT=2\pi$, so
that the long-time evolution of the velocity would also be close to Brownian. This is
what we shall show.

An intimately related issue is the relative motion of several particles, released in the
same realization of the wave field. Even though each particle velocity diffuses for $t
\in [0,2\pi]$, their motions are not independent. We shall also show that for large $A$
the motions of any finite family of particles released at initial data $(Q_0^{(\nu)},
P_0^{(\nu)})$, $1 \le \nu \le \cN$, approaches a family of $\cN$ independent processes.
This can also be expected from the consideration of the top Lyapunov exponent of the
dynamics \eqref{eqq2b}--\eqref{eqp2b} in the limit $A \to \infty$.

\begin{theorem} \label{NpartGauss}
  For any $\cN>0$, the momentum processes $P^{(\nu)}$ defined by
  \begin{eqnarray}
    \rmd Q_t^{(\nu)}
    &=& A P_t^{(\nu)} \ \rmd t \, , \  \
    Q_0^{(\nu)} = q_0^{(\nu)} \, ,
    \\
    \rmd P_t^{(\nu)}
    &=& (\sin Q_t^{(\nu)}) \rmd C_t + (\cos Q_t^{(\nu)}) \rmd S_t \, ,\  \
    P_0^{(\nu)} = p_0^{(\nu)} \, ,
  \end{eqnarray}
  with $\cN$ different initial data
  $(q_0^{(\nu)}, p_0^{(\nu)}) \in \TT \times \RR$, $1 \le \nu \le \cN$,
  converge as $A \to \infty$ to $\cN$ independent Wiener processes
  with variance $\pi t$,
  and convergence is in law in $C(\RR_+,\RR^\cN)$.
\end{theorem}

The key argument in the proof is the following weak convergence theorem, where we write
now $n = \pi^{1/2} A$. Consider the two--dimensional diffusion process indexed by $n \ge
1$, solution of the SDE on $\RR^2$
\begin{equation} \label{eds}
  \left\{
  \begin{aligned}
    \frac{\rmd U^n_t}{\rmd t}
    &= n V^n_t \, ,\
    U_0 = u \, ,
  \\
    \rmd V^n_t
    &= \sin(U^n_t) \rmd W_t \, ,\
    V_0 = v \, ,
  \end{aligned}
  \right.
\end{equation}
where $(u,v)\not\in\{(k\pi,0),\ k\in\ZZ\}$ and $W$ is a standard Brownian motion. We
prove
\begin{theorem}\label{theth}
  As $n\to\infty$,
    $$V^n \Rightarrow v + \frac{1}{\sqrt{2}} B \, ,$$
  where $\{B_t,\ t \ge 0\}$ is a standard one--dimensional Brownian motion,
  and the convergence is in law in $C(\RR_+,\RR)$.
\end{theorem}

\section{A change of time scale}

Note that for any $n  \ge  1$, the law of $\{(U^n_t,V^n_t),\ t  \ge  0 \}$, the solution
of \eqref{eds}, is characterized by the statement
\begin{equation*}
\left\{
\begin{aligned}
  \frac{\rmd U^n_t}{\rmd t}
  &=
  nV^n_t
  ,\ U_0 = u \, ,
  \\
  V^n
  &\text{ is a martingale},\
  \frac{\rmd \langle V^n \rangle_t}{\rmd t}
  =
  \sin^2(U^n_t)
  \, ,\
  V^n_0 = v  \, .
\end{aligned}
\right.
\end{equation*}
Now define (like B\'enisti and Escande \cite{BE97,BE98})
\begin{equation}
  X_t=U^n_{n^{-2/3}t} \, ,\quad Y_t=n^{1/3}V^n_{n^{-2/3}t} \, .
  \label{XYscale}
\end{equation}
We first note that $X_0=u$, $Y_0=n^{1/3}v$, $Y$ is a martingale, and
\begin{equation*}
\left\{
\begin{aligned}
  \frac{\rmd X_t}{\rmd t}
  &=
  n^{-2/3} \frac{\rmd U^n}{\rmd t} (n^{-2/3}t)
  =
  n^{1/3}V^n_{n^{-2/3}t}=Y_t
  \, , \\
  \langle Y \rangle_t
  &=
  n^{2/3} \langle V^n \rangle_{n^{-2/3}t}
  \, ,\
  \frac{\rmd \langle Y \rangle_t}{\rmd t}
  = \sin^2(X_t) \, .
\end{aligned}
\right.
\end{equation*}
Using a well--known martingale representation theorem, we can pretend that there exists
a standard Brownian motion $\{B_t,\ t \ge  0\}$ such that
\begin{equation}\label{equ}
\left\{
\begin{aligned}
  \frac{\rmd X_t}{\rmd t}
  &= Y_t \, ,\
  X_0=u \, ,\\
  \rmd Y_t
  &= \sin(X_t) \rmd B_t \, ,\
  Y_0=n^{1/3}v \, .
\end{aligned}
\right.
\end{equation}
Note that the process $\{(X_t,Y_t),\ t \ge 0\}$ still depends upon $n$, but only through
the value of $Y_0$.

On the other hand, $V^n_t=n^{-1/3}Y_{n^{2/3}t}$. Hence
  $$V^n_t = v + n^{-1/3} \int_0^{n^{2/3}t} \sin(X_s) \rmd B_s \, ,$$
in other words $V^n$ is a martingale such that $V^n_0=v$ and
  $$\langle V^n \rangle_t = n^{-2/3} \int_0^{n^{2/3}t} \sin^2(X^n_s) \rmd s \, .$$
Here we recall the fact that the process $X$ depends upon $n$ (through the initial
condition of $Y$), unless $v=0$. Consequently
\begin{equation}
\label{limit_n}
  \lim_{n\to\infty}\langle V^n\rangle_t
  =
  t \times \lim_{n\to\infty}\frac{1}{n^{2/3}t}\int_0^{n^{2/3}t}\sin^2(X^n_s)\rmd s
\end{equation}
and in order to prove Theorem \ref{theth} it suffices to show that the above limit is
$t/2$.

\section{Qualitative properties of the solution of \eqref{equ}}

We now consider the two--dimensional diffusion process
\begin{equation}\label{syst-diff-stoch}
\left\{
\begin{aligned}
  \frac{\rmd X_t}{\rmd t}
  &=Y_t \, ,\
  X_0=x \, ,
  \\
  \rmd Y_t
  &=\sin(X_t) \rmd B_t \, ,\
  Y_0=y \, ,
\end{aligned}
\right.
\end{equation}
with values in the state--space $E=[0,2\pi)\times\RR\backslash \{(0,0),(\pi,0)\}$, where
$2\pi$ is identified with 0. We first prove that the process $\{(X_t,Y_t),\ t \ge 0\}$
is a conservative $E$--valued diffusion. Indeed,
\begin{proposition}
Whenever the initial condition $(x,y)$ belongs to $E$, $$\inf\{t>0,\ (X_t,Y_t)\in
\{(0,0),(\pi,0)\}\}=+\infty\quad\text{a.s.}$$
\end{proposition}
\bpf We define the stopping time
  $$\tau=\inf\{t,\ (X_t,Y_t)=(0,0)\} \, .$$
Let $R_t=X^2_t+Y^2_t$, $Z_t=\log R_t$, $t \ge 0$. A priori, $Z_t$ takes its values in
$[-\infty,+\infty)$. It\^o calculus on the interval $[0,\tau)$ yields
\begin{eqnarray*}
  \rmd X^2_t
  &=& 2 X_t Y_t \rmd t \, ,
  \\
  \rmd Y^2_t
  &=& 2 Y_t \sin(X_t) \rmd B_t + \sin^2(X_t) \rmd t \, ,
  \\
  \rmd Z_t
  &=& \frac{\rmd R_t}{R_t}-\frac{\rmd \langle R \rangle _t}{2 R^2_t} \\
  &=& \frac{2Y_tX_t+\sin^2(X_t)}{R_t}\rmd t
     - 2 \frac{Y^2_t \sin^2(X_t)}{R^2_t}\rmd t 
     + 2 \frac{Y_t\sin(X_t)}{R_t} \rmd B_t \, .
\end{eqnarray*}
Now clearly $|\sin(x)|\le |x|$, $\sin^2(x)\le x^2$, and it follows from the above and
standard inequalities that on the time interval  $[0,\tau)$,
  $$Z_t \ge  Z_0-2t+\int_0^t \varphi_s \rmd B_s \, ,$$
where $|\varphi_s|\le 1$. Hence the process $\{Z_t,\ t \ge 0\}$ is bounded from below on
any finite time interval, which implies that $\tau=+\infty$ a.s., since $\tau=\inf\{t,
Z_t =-\infty\}$. A similar argument shows that $\tau'=+\infty$ a.s., where
  $$\tau'=\inf\{t,\ (X_t,Y_t)\in\{(0,0),(\pi,0)\}\} \, .$$ \epf

We next prove (here and below $\cB_E$ stands for the $\sigma$--algebra of Borel subsets
of $E$) the
\begin{proposition}\label{malliavin}
The transition probabilities
  $$ \{p((x,y);t,A)
     := \PP_{x,y}((X_t,Y_t)\in A),\ (x,y)\in E,\ t>0,\ A\in\cB_E\}$$
have smooth densities $p((x,y);t,(x',y'))$ with respect to Lebesgue's measure $\rmd x'
\rmd y'$ on $E$.
\end{proposition}
\bpf Consider the Lie algebra of vector fields on $E$ generated by
$X_1=\sin(x)\frac{\partial}{\partial y}$, $X_2=[X_0,X_1]$ and $X_3=[X_0,[X_0,X_1]]$,
where $X_0=y\frac{\partial}{\partial x}$. This Lie algebra has rank 2 at each point of
$E$. The result is now a  standard consequence of the well--known Malliavin calculus,
see e.g.\ Nualart \cite{nu}. \epf

 \begin{proposition}\label{irreduct}
The $E$--valued diffusion process $\{(X_t,Y_t),\ t \ge 0\}$ is topologically
irreducible, in the sense that for all $(x,y)\in E$, $t>0$, $A\in\cB_E$ with non empty
interior,
  $$\PP_{x,y}((X_t,Y_t)\in A)>0 \, .$$
\end{proposition}
\bpf From Stroock--Varadhan's support theorem, see e.g.\ Ikeda--Watanabe \cite{iw}, the
support of the law of $(X_t,Y_t)$ starting from $(X_0,Y_0)=(x,y)$ is the closure of the
set of points which the following controlled ordinary differential equation can reach at
time $t$ by varying the control $\{u(s),\ 0\le s\le t\}$ in the class of piecewise
continuous functions~:
\begin{equation}
\left\{
\begin{aligned}\label{cont-ode}
  \frac{\rmd x}{\rmd s}(s)&=y(s) \, ,\quad x(0)=x
  \, ;\\
  \frac{\rmd y}{\rmd s}(s)&=\sin(x(s)) u(s) \, ,\quad y(0)=y
  \, .
\end{aligned}
\right.
\end{equation}
It is not hard to show that the set of accessible points at time $t>0$ by the solution
of \eqref{cont-ode} is dense in $E$. The result now follows from the fact that the
transition probability is absolutely continuous with respect to Lebesgue's measure, see
Proposition \ref{malliavin}. \epf



We next prove the
\begin{lemma}
$$\PP\left(|Y_t|\to\infty,\ \mathrm{as }\  t\to\infty\right)=0 \, .$$
\end{lemma}
\bpf
 The Lemma follows readily from the fact that
$$Y_t=W\left(\int_0^t\sin^2(X_s)\rmd s\right),$$ where $\{W(t),\ t \ge 0\}$ is a scalar
Brownian motion. Then either $\int_0^t\sin^2(X_s)\rmd s$ is bounded and $Y_t$ is
finite, or the integral diverges and $Y_t$ is finite anyway because $W$ is recurrent.
\epf

Hence the topologically irreducible $E$--valued Feller process $\{(X_t,Y_t), \ t \ge
0\}$ is recurrent. Its unique (up to a multiplicative constant) invariant measure is the
Lebesgue measure on $E$, so that in particular the process is null--recurrent.  It then
follows from (ii) in Theorem 20.21 from Kallenberg \cite{ka}
\begin{lemma}\label{petit}
For all $M>0$, as $t\to\infty$,
  $$\frac{1}{t} \int_0^t{\bf1}_{\{|Y_s| \le M\}}\rmd s \to 0 \quad \text{a.s.}$$
\end{lemma}

\section{A path decomposition of the process $\{(X_t,Y_t),\ t \ge 0\}$}

We first define two sequences of stopping times. Let $T_0=0$ and
\begin{align*}
\text{for }\ell\text{ odd},& \quad T_\ell=\inf\{t>T_{\ell-1},\ |Y_t| \ge  M+1\}
  \, ,\\
\text{for }\ell\text{ even},& \quad T_\ell=\inf\{t>T_{\ell-1},\ |Y_t|\le M\} \, .
\end{align*}
Let now $\tau_0=T_1$. We next define recursively $\{\tau_k,\ k \ge 1\}$ as follows.
Given $\tau_{k-1}$, we first define $$L_k=\sup\{\ell \ge 0,\ \tau_{k-1} \ge
T_{2\ell+1}\} \, .$$ Now let
  $$ \eta_k
     = \begin{cases}
         \tau_{k-1} &\text{if $\tau_{k-1} <    T_{2L_k+2}$} \, ,  \\
         T_{2L_k+3} &\text{if $\tau_{k-1} \ge  T_{2L_k+2}$} \, .
       \end{cases}$$
We now define
  $$ \tau_k = \inf\{t>\eta_k,\ |X_t-X_{\eta_k}|=2\pi\}
              \wedge \inf\{t>\eta_k,\ |Y_t-Y_{\eta_k}|>1\} \, .$$
It follows from the above definitions that
  $$ \int_0^t {\bf1}_{\{|Y_s| \ge  M+1\}}\sin^2(X_s)\rmd s
     \le \sum_{k=1}^\infty\int_{\eta_k\wedge t}^{\tau_k\wedge t}\sin^2(X_s)\rmd s
     \le \int_0^t\sin^2(X_s)\rmd s \, , $$
a statement which will be refined in the proof of Proposition \ref{prop_final}. Define
\begin{align*}
  K^0 &= \{k \ge 1,\ |Y_{\tau_k}-Y_{\eta_k}|<1\} \, , \\
  K^1 &= \{k \ge 1,\ |Y_{\tau_k}-Y_{\eta_k}|=1\} \, , \\
  K^{\ }_t &= \{k \ge 1,\ \eta_k< t\} \, , \\
  K^0_t &= K^0\cap K_t \, , \\
  K^1_t &= K^1\cap K_t \, .
\end{align*}
We first prove the
\begin{lemma}\label{neglig}
$$\frac{1}{t}\sum_{k\in K^1_t}(\tau_k-\eta_k)\to 0$$ in $L^1(\Omega)$ as $M\to\infty$,
uniformly in $t>0$.
\end{lemma}
\bpf We shall use repeatedly the fact that  since $|Y_{\eta_k}| \ge  M>2$,
$|Y_{\eta_k}|- 1 \ge  |Y_{\eta_k}|/2$. We have that (see the Appendix below), since
$\tau_k-\eta_k\le4\pi/|Y_{\eta_k}|$,
\begin{eqnarray*}
  \PP(k\in K^1|\cF_{\eta_k})
  & \le &
  \PP\left(\sup_{\eta_k \le t \le \tau_k} |Y_t-Y_{\eta_k}| \ge 1
          \Big| \cF_{\eta_k} \right)
  \\
  & \le &
  2 \exp(-|Y_{\eta_k}|/(8\pi)) \, .
\end{eqnarray*}
Consequently, using again the inequality $\tau_k-\eta_k\le4\pi/|Y_{\eta_k}|$, we deduce
that
\begin{eqnarray*}
  \EE\left[(\tau_k-\eta_k){\bf1}_{\{k \in K^1\}}|\cF_{\eta_k}\right]
  & \le &
  \frac{8\pi}{|Y_{\eta_k}|}\exp(-|Y_{\eta_k}|/(8\pi))
  \\
  & \le &
  \frac{8\pi}{|Y_{\eta_k}|}\exp(-M/(8\pi)) \, .
\end{eqnarray*}
On the other hand, whenever $k\in K^0$,
  $$\tau_k-\eta_k \ge 2\pi/(|Y_{\eta_k}|+1) \ge \pi/|Y_{\eta_k}| \, .$$
Now, provided $t \ge  4\pi/M$,
\begin{align*}
  2t
  & \ge
  t+\frac{4\pi}{M}
  \\ &
  \ge \EE\left[\sum_{k\in K^0_t}(\tau_k-\eta_k)\right]
  \\ &
  \ge \pi\EE\left[\sum_{k\in K_t}{\bf1}_{\{k\in K^0\}}\frac{1}{|Y_{\eta_k}|}\right]
  \\ &
  \ge \frac{\pi}{2}\EE\left[\sum_{k\in K_t}\frac{1}{|Y_{\eta_k}|}\right] \, ,
\end{align*}
since
\begin{align*}
\PP(k\in K^0|\cF_{\eta_k})&=1-\PP(k\in K^1|\cF_{\eta_k})\\ & \ge 1-2\exp(-M/(8\pi))\\ &
\ge 1/2 \, ,
\end{align*}
provided $M$ is large enough. Finally
\begin{align*}
\frac{1}{t}\EE\left[\sum_{k\in K^1_t}(\tau_k-\eta_k)\right]&\le 32\exp(-M/(8\pi))
\frac{\EE\left[\sum_{k\in K_t}|Y_{\eta_k}|^{-1}\right]}{\EE\left[\sum_{k\in
K_t}|Y_{\eta_k}|^{-1}\right]}\\ &=32\exp(-M/(8\pi))\\ &\to 0 \, ,
\end{align*}
as $M\to\infty$, uniformly in $t$. \epf

Now, for any $k\in K^0$,
\begin{eqnarray}
  \int_{\eta_k}^{\tau_k} \sin^2(X_s) \rmd s
  &=&
  \frac {\tau_k-\eta_k} {2\pi} \int_0^{2\pi} \sin^2(x) \rmd x
  \nonumber \\
  && + \int_{\eta_k}^{\tau_k} \sin^2(X_s)
                \left[1-\frac{Y_s(\tau_k-\eta_k)}{2\pi}\right]
                \rmd s \, ,
  \label{intsin2}
\end{eqnarray}
and we have
\begin{align*}
\left|\int_{\eta_k}^{\tau_k}\sin^2(X_s)\left[1-\frac{Y_s(\tau_k-\eta_k)}{2\pi}\right]\rmd
s \right| &=\left|\int_{\eta_k}^{\tau_k}\int_{\eta_k}^{\tau_k}\sin^2(X_s) \frac{Y_r -
Y_s}{2\pi}  \rmd r \rmd s \right| \\
&\le\frac{1}{2\pi}\int_{\eta_k}^{\tau_k}\int_{\eta_k}^{\tau_k}|Y_r-Y_s|\rmd r\rmd s \, .
\end{align*}

Finally we have the
\begin{lemma}\label{errorterm}
  Uniformly in $t>0$,
  $$ \frac{\sum_{k\in K^0_t} \int_{\eta_k}^{\tau_k} \int_{\eta_k}^{\tau_k}
                             |Y_r-Y_s| \rmd r \rmd s}
          {\sum_{k \in K^0_t} (\tau_k-\eta_k)}
     \to 0 $$
  a.s., as $M\to\infty$.
\end{lemma}
\bpf Since $|Y_t-Y_{\eta_{k}}|\le1$ for $\eta_{k}\le t\le\tau_k$,
\begin{align*}
  \frac{\sum_{k\in K^0_t} \int_{\eta_k}^{\tau_k}\int_{\eta_k}^{\tau_k}
                               |Y_r-Y_s| \rmd r \rmd s}
       {\sum_{k\in K^0_t}(\tau_k-\eta_k)}
  & \le 2 \sup_{k\in K^0_t}(\tau_k-\eta_{k})
  \\
  & \le 8 \pi/M
  \\
  & \to 0 \, ,
\end{align*}
as $M\to\infty$, uniformly in $t$. \epf

We are now in a position to prove the following ergodic type theorem, from which Theorem
\ref{theth} will follow~:
  \begin{proposition}\label{prop_final}
As $t\to\infty$, $$\frac{1}{t}\int_0^t\sin^2(X_s)\rmd s\to\frac{1}{2}$$ in probability.
\end{proposition}
\bpf We first note that $$[0,t]=B^0_t\cup B^1_t\cup C_t \, ,$$ where
\begin{align*}
  B^0_t   &= [0,t]\cap\left(\cup_{k\in K^0_t}[\eta_k,\tau_k]\right) \, ,\\
  B^1_t   &= [0,t]\cap\left(\cup_{k\in K^1_t}[\eta_k,\tau_k]\right) \, ,\\
  C^{ }_t &= [0,t]\backslash(B^0_t\cup B^1_t) \, .
\end{align*}
We have
\begin{eqnarray*}
  \frac{1}{t} \int_0^t \sin^2(X_s)\rmd s
  & = &
  \frac{1}{t} \int_0^t {\bf1}_{B^0_t}(s) \sin^2(X_s) \rmd s
  + \frac{1}{t} \int_0^t {\bf1}_{B^1_t}(s) \sin^2(X_s) \rmd s
  \\ & &
  + \frac{1}{t} \int_0^t {\bf1}_{C_t}(s) \sin^2(X_s) \rmd s \, .
\end{eqnarray*}

Now $C_t\subset\{s\in[0,t],\ |Y_s|\le M+1\}$, so for each fixed $M>0$, it follows from
Lemma \ref{petit} that the  last term can be made arbitrarily small, by choosing $t$
large enough.  The second term goes to zero as $M\to\infty$, uniformly in $t$, from
Lemma \ref{neglig}. Finally the first term equals the searched limit, plus an error term
which goes to 0 as $M\to\infty$, uniformly in $t$, see Lemma \ref{errorterm} and the
following fact, which follows from the combination of Lemma \ref{neglig} and Lemma
\ref{petit} : $$ \frac{1}{t}\sum_{k\in K^0_t}(\tau_k-\eta_k)\to 1$$ in probability, as
$n\to\infty$. \epf

We can finally proceed with the

\noindent{\sc Proof of Theorem \ref{theth}} All we have to show is that (see
\eqref{limit_n})
  $$\lim_{n\to\infty}\frac{1}{n^{2/3}t}\int_0^{n^{2/3}t}\sin^2(X^n_s)\rmd s
      = \frac{1}{2\pi} \int_0^{2\pi}\sin^2(x')\rmd x'
      = \frac{1}{2}$$
in probability. In the case
$v=0$, the process $\{(X^n_t,Y^n_t)\}$ does not depend upon $n$, and the result follows
precisely from Proposition \ref{prop_final}. Now suppose that $v\not=0$. In that case,
the result can be reformulated equivalently as follows. For some $x\in\RR$, $y\not=0$,
each $t>0$, define the process $\{(X^t_s,Y^t_s),\ 0\le s\le t\}$ as the solution of the
SDE
\begin{equation*}
\left\{
\begin{aligned}
  \frac{\rmd X^t_s}{\rmd s} &= Y^t_s \, ,\ X^t_0=x \, ,
  \\ \rmd Y^t_s &= \sin(X^t_s) \, \rmd W_s \, ,\
  Y^t_0 = \sqrt{t} \, y \, .
\end{aligned}
\right.
\end{equation*}
We need to show that
  $$\frac{1}{t}\int_0^t\sin^2(X^t_s)\rmd s \to
    \frac{1}{2\pi}\int_0^{2\pi}\sin^2(x')\rmd x'$$ in probability, as $t\to\infty$. Note that
in time $t$, the process $Y^t$ starting from $\sqrt{t}y$ can come back near the origin.

It is easily seen, by introducing the Markov time $\tau^t_M=\inf\{s>0,\ |Y^t_s|\le M\}$
and exploiting the strong Markov property, that
  $$\frac{1}{t}\int_0^t{\bf1}_{\{|Y^t_s|\le M\}}\rmd s\to0\quad \text{a.s.}$$
follows readily from Lemma \ref{petit}. The rest of the argument leading to Proposition
\ref{prop_final} is based upon limits as $M\to\infty$, uniformly with respect to $t$. It
thus remains to check that the fact that $Y^t_0$ now depends upon $t$ does not spoil
this uniformity, which is rather obvious. \epf

\begin{remark}{\rm{This proof holds uniformly with respect to initial data $(u,v)$ satisfying $|v| \geq a$ for any $a>0$.
}}\end{remark}

\section{Proof of Theorem \ref{NpartGauss}}
\label{proofNpart}

Our Theorem \ref{NpartGauss} now appears as a simple corollary of Theorem \ref{theth}.

\bpf We first prove the Theorem for $t \in [0,2\pi]$. Then the vector ${\mathbf{P}} =
(P^{(1)}, \ldots, P^{(\cN)})$ is a martingale in $\RR^\cN$, and to prove our claim it
suffices to show that its quadratic variation matrix converges to $\pi t$ times the
identity matrix. The diagonal elements of the matrix are
\begin{equation*}
  \qvar{P^{(\nu)}}_t
    = \int_0^t (\sin^2 Q^{(\nu)}_s + \cos^2 Q^{(\nu)}_s) \pi \rmd s = \pi t
\end{equation*}
and we only need to compute the cross-variation
\begin{eqnarray}
  \qvar{P^{(\nu)}, P^{(\nu')}}_t
  &=& \int_0^t (\sin Q^{(\nu)}_s \sin Q^{(\nu')}_s + \cos Q^{(\nu)}_s \cos Q^{(\nu')}_s) \pi \rmd s
  \nonumber
  \\
  &=& \int_0^t \cos(Q^{(\nu)}_s - Q^{(\nu')}_s) \pi \rmd s \, .
\end{eqnarray}
Now, define (with $n = \pi^{1/2} A$)
\begin{eqnarray*}
    U^n_t = {\frac 1 2} (Q^{(\nu)}_t - Q^{(\nu')}_t) & , &
    V^n_t = n^{-1} \frac{\rmd U^n_t}{\rmd t} \, ,
    \\
    U'^n_t = {\frac 1 2} (Q^{(\nu)}_t + Q^{(\nu')}_t) & , &
    V'^n_t = n^{-1} \frac{\rmd U'^n_t}{\rmd t} \, .
\end{eqnarray*}
These processes solve the stochastic differential equation
\begin{eqnarray}
  \rmd U^n_t &=& n V^n_t \, \rmd t \,  ,\   \quad
  U^n_0 =  \frac {q^{(\nu)}_0 - q^{(\nu')}_0} 2 \, ,
  \label{UJ1}
  \\
  V^n_0 &=&  \frac {p^{(\nu)}_0 - p^{(\nu')}_0} {2 \pi^{1/2}} \, ,
  \label{UJ2}
  \\
  \rmd V^n_t
  &=& \frac {1}{2\sqrt{\pi}}
      \left(\sin(U'^n_t + U^n_t) - \sin(U'^n_t - U^n_t)\right) \rmd C_t
  \nonumber \\
  && + \frac {1}{2\sqrt{\pi}}
       \left(\cos(U'^n_t + U^n_t) - \cos(U'^n_t - U^n_t)\right) \rmd S_t
  \nonumber \\
  &=& \sin U^n_t \ \rmd W^n_t
  \label{UJ3}
\end{eqnarray}
where the process $W^n$ is the martingale defined by $W^n_0=0$ and
  $$ \rmd W^n_t = \pi^{-1/2} (\cos U'^n_t \ \rmd C_t - \sin U'^n_t \ \rmd S_t)\, . $$
The quadratic variation of $W^n$ is
  $$ \qvar{W^n}_t = \pi^{-1} \int_0^t (\cos^2 U'^n_s + \sin^2 U'^n_s)\ \pi \ \rmd s
                = t$$
in view of the quadratic and cross variations \eqref{EECS} of $(C,S)$, and this result
does not depend on the process $U'^n$ (which follows the center of mass of the two
particles $\nu$ and $\nu'$). Thus $W^n$ is a standard Wiener process, and $(U^n,V^n)$
defined by \eqref{UJ1}--\eqref{UJ2}--\eqref{UJ3} satisfies the hypotheses of Theorem
\ref{theth}. Hence, with $X$ defined by \eqref{XYscale},
\begin{eqnarray*}
  \qvar{P^{(\nu)}, P^{(\nu')}}_t
  &=& \int_0^t \cos(2 U^n_s) \pi \rmd s
  = \int_0^t (1 - 2 \sin^2 U^n_s) \pi \rmd s
  \\
  &=& \pi t \left( 1 - \frac 2 { n^{2/3} t} \int_0^{n^{2/3} t} \sin^2 X^n_{s'} \rmd s'
            \right)
\end{eqnarray*}
which converges in probability to 0 for $n \to \infty$ as shown in the proof of the
theorem.

Now we consider the process over the interval $[0,4\pi]$, taking into account that
$(C,S)$ over the whole interval is neither a martingale nor Markov. From the given
initial data, amplitude $A$ and realization of $(C,S)$, we define a subsidiary set of
particles $\cN+1 \leq \nu \leq 2\cN$, with initial data
  $$Q^{(\nu)}_0 = Q^{(\nu-\cN)}_{2\pi} \, ,\ P^{(\nu)}_0 = P^{(\nu-\cN)}_{2\pi} \, ,$$
of which a.s.\ none coincides (modulo $2\pi$ for $q$) with any of the initial data
$(Q^{(\nu)}_0, P^{(\nu)}_0)$. Recalling that at $t=2\pi$ the law of any
$P^{(\nu)}_{2\pi}$ (for $1 \le \nu \le \cN$) is Gaussian with variance $2 \pi^2$ (so
that its probability density is bounded by $1/\sqrt{2 \pi \ 2 \pi^2} < 1/(2\pi)$), this
ensures that, given $\epsilon
> 0$,
\begin{eqnarray*}
  && \PP \left( \min_{1 \le \nu, \nu' \le \cN} | P^{(\nu)}_0 - P^{(\nu')}_{2\pi} | < \epsilon
                      \  {\textrm{or}} \
                       \min_{1 \le \nu < \nu' \le \cN} | P^{(\nu)}_{2\pi} - P^{(\nu')}_{2\pi} | < \epsilon
              \right)
    \\
    & \le & \frac {3 \epsilon \cN^2} {\sqrt{2 \pi \ 2 \pi^2}}
    < \epsilon \cN^2 / 2
  \,  .
\end{eqnarray*}
For processes such that $\min_{\nu,\nu'} | P^{(\nu)}_0 - P^{(\nu')}_{2\pi} | \ge
\epsilon$ and $\min_{\nu<\nu'} | P^{(\nu)}_{2\pi} - P^{(\nu')}_{2\pi} | \ge \epsilon$,
the above proof holds (uniformly with respect to initial data) for the full set of
$2\cN$ particles defined here. It follows that a.s.\ the $P^{(\nu)}_t - P^{(\nu)}_0$, $t
\in [0,2\pi]$, $1 \leq \nu \leq 2\cN$ are mutually independent Wiener processes, each
with variance $\pi t$. This implies that $P^{(\nu)}_t - P^{(\nu)}_0$, $t \in [0,4\pi]$,
$1 \leq \nu \leq \cN$ are also mutually independent Wiener processes.

For any interval $[0,2k\pi]$ with $k \in \NN_0$ the same argument reduces the
$\cN$-particle problem to $k\cN$ particles over $[0,2\pi]$, and the proof is complete.
\epf

\begin{remark}
\rm{Our theorem allows that $P^{(\nu)}_0 = P^{(\nu')}_0$
  for some $\nu \neq \nu'$ with $1 \le \nu, \nu' \le \cN$.}
\end{remark}

\section{Perspectives}
\label{persp}

The implications of Theorem \ref{NpartGauss} are twofold. First, for $\cN=1$, they
support the observation \cite{BE97,BE98,CEV,E08} that, in systems with finite $\cM \gg
1$ and $A_0/{\mfm} \gg 1$ (with appropriate scaling), the long-time behaviour of a
single particle in a periodic wave field exhibits statistical properties approaching
those of the Brownian motion. To complete the connection with physics literature, one
must now discuss how the finite sums in \eqref{dotv} approach the right hand side of
\eqref{eqp2b}, and how this implies that the solutions of the first equation approach
the solutions of the latter equation \cite{su}. This will be discussed separately.

Second, and this is conceptually more fundamental, for a single sample of the Wiener
wavefield (with $\cM = \infty$ formally), with strong spatial correlations (thanks to
the single wavevector $k_0$ in the model), the limit $A_0/{\mfm} \to \infty$ leads to
independence of the evolutions of all $\cN$ particles, which can then be collectively
described by the diffusion equation. While such an independence is often admitted
without proof in physics practice, our work provides an explicit justification to it. We
even prove a little more than usual one-- or two--time statements, as in our limit
$p^{(\cN)}$ is independent in law from even the full evolution data $\{p^{(\nu)}(t), \ 0
\leq t \leq T, \ 1 \leq \nu \leq \cN-1\}$.

This second implication is an important issue, as the acceleration of a passive particle
is a Hamiltonian process while the diffusion equation is irreversible. While the key to
this irreversibility is clearly the fact that the diffusion process only relates to the
momentum component of particle evolution, we shall further investigate the interplay of
limits $\cN \to \infty$ and $A \to \infty$, or $\mfm \to 0$, with stochasticity versus
Hamiltonian ``conservativeness'' in $(Q,P)$ variables in future work.

Finally, we followed general practice in discussing the stochastic acceleration problem
in only one space dimension \cite{Chirikov79,HW04}. This is rather classical, and it
applies \eg to particle motion along magnetic field lines in strongly magnetized
plasmas~; higher-dimensional motions may call for different elementary models.


\appendix
\section{}
\label{app}

For the convenience of the reader, we prove the following well--known result (see e.g.\
ex.\ IV.3.16 in \cite{RevuzYor99})
\begin{proposition}
Let $\eta$ and $\tau$ be two stopping times such that $0\le\eta\le\tau\le\eta+T$ and
$M_t=\int_0^t\varphi_s \rmd B_s$, where $\{B_t,\ t \ge 0\}$ is a standard Brownian
motion and $\{\varphi_t,\ t \ge 0\}$ is progressively measurable and satisfies $\int_0^t
\varphi_s^2 \rmd s \le k^2 t$ for all $t \ge 0$. Then for all $b>0$,
  $$\PP\left(\sup_{\eta\le t\le\tau} |M_t-M_{\eta}| \ge  b\right)
     \le2 \exp\left(-\frac{b^2}{2k^2T}\right).$$
\end{proposition}
\bpf From the optional stopping theorem, it suffices to treat the case $\eta=0$,
$\tau=T$. We have
  $$ \PP(\sup_{0\le t\le T}|M_t| \ge  b)
     \leq
     \PP(\sup_{0\le t\le T}M_t \ge b)
       + \PP(\inf_{0\le t\le T}M_t\le -b) \, . $$
We estimate the first term on the right. The second one is bounded by the same quantity.
Define for all $\lambda>0$
  $$\cM^\lambda_t
    = \exp\left(\lambda M_t-\frac{\lambda^2}{2}\int_0^t\varphi^2_s\rmd s\right).$$
Then
\begin{align*}
  \PP(\sup_{0\le t\le T}M_t \ge  b)
  &\le \PP\left(\sup_{0\le t\le T}\cM^\lambda_t
                \ge \exp(\lambda b-\lambda^2k^2T/2)\right)
  \\
  &\le \exp(\lambda^2k^2T/2-\lambda b) \, ,
\end{align*}
from Doob's inequality, since $\{\cM^\lambda_t,\ t \ge 0\}$ is a martingale with mean
one. Optimizing the value of $\lambda$, we deduce that
  $$\PP \left( \sup_{0\le t\le T}M_t \ge b \right)
    \le \exp \left(-\frac{b^2}{2k^2T} \right),$$
from which the result follows. \epf

\section*{Acknowledgements}

This work stems from many stimulating discussions with Dominique Escande, Fabrice Doveil
and members of the Turbulence plasma team. Comments from Claude Bardos, Nils Berglund
and C\'edric Villani were useful. Jonathan Mattingly was helpful in providing the
reference for Lemma \ref{petit}. YE was partly supported by CNRS through a delegation
position.


\end{document}